\documentclass[11pt, reqno]{amsart}
\usepackage{amssymb}
\usepackage{amscd}
\usepackage{verbatim,ifthen}
\usepackage{color}
\usepackage{latexsym}
\usepackage{tikz}
\usepackage{tikz-cd}
\usepackage{mathrsfs}
\usepackage{wrapfig}
\usetikzlibrary{shapes}
\usepackage{color}
\usetikzlibrary{arrows.meta}
\usepackage{bbm}
\usetikzlibrary{matrix}
\usetikzlibrary{calc}
\usetikzlibrary{arrows,intersections}
\usepackage{pgfplots}
\usepackage{multicol}
\usepackage{array}
\newcolumntype{M}[1]{>{\centering\arraybackslash}m{#1}}

\usepackage[colorlinks, linkcolor=black, citecolor=magenta, linktocpage]{hyperref}
\addtocontents{toc}{\setcounter{tocdepth}{1}}

\usepackage{amsmath}

\numberwithin{equation}{section}

\let\oldtocsection=\tocsection
 
\let\oldtocsubsection=\tocsubsection

\renewcommand{\tocsection}[2]{\hspace{0em}\oldtocsection{#1}{#2}}
\renewcommand{\tocsubsection}[2]{\hspace{1em}\oldtocsubsection{#1}{#2}}

\def\XXint#1#2#3{{\setbox0=\hbox{$#1{#2#3}{\int}$ }
\vcenter{\hbox{$#2#3$ }}\kern-.6\wd0}}

\usepackage{eucal}
\usepackage{calc}  
\usepackage{enumitem} 
\usepackage{tensor}
\usepackage{graphicx,wrapfig,lipsum}
\usepackage{etoolbox}
\usepackage{marginnote}
\usepackage{lipsum}
\makeatletter
\patchcmd{\@mn@margintest}{\@tempswafalse}{\@tempswatrue}{}{}
\patchcmd{\@mn@margintest}{\@tempswafalse}{\@tempswatrue}{}{}
\reversemarginpar 
\makeatother
\usepackage{scrextend}

\makeatletter
\DeclareRobustCommand\widecheck[1]{{\mathpalette\@widecheck{#1}}}
\def\@widecheck#1#2{%
    \setbox\z@\hbox{\m@th$#1#2$}%
    \setbox\tw@\hbox{\m@th$#1%
       \widehat{%
          \vrule\@width\z@\@height\ht\z@
          \vrule\@height\z@\@width\wd\z@}$}%
    \dp\tw@-\ht\z@
    \@tempdima\ht\z@ \advance\@tempdima2\ht\tw@ \divide\@tempdima\thr@@
    \setbox\tw@\hbox{%
       \raise\@tempdima\hbox{\scalebox{1}[-1]{\lower\@tempdima\box
\tw@}}}%
    {\ooalign{\box\tw@ \cr \box\z@}}}
\makeatother
\title{Remarks on the Quadratic Orthogonal Bisectional Curvature}
\author{Kyle Broder}
\address[Kyle Broder]{Mathematical Sciences Institute, Australian National University, Acton, ACT 2601, Australia; BICMR, Peking University, Beijing, 100871, People's Republic of China}
\email{kyle.broder{@}anu.edu.au}
\thanks{Supported by an Australian Government Research Training Program Scholarship}
\begin{document}

\maketitle

\begin{abstract}
We exhibit a curious link between the Quadratic Orthogonal Bisectional Curvature, combinatorics, and distance geometry.  The Weitzenb\"ock curvature operator, acting on real $(1,1)$--forms, is realized as the Dirichlet energy of a finite graph, weighted by a matrix of the curvature. These results also illuminate the difference in the nature of the Quadratic Orthogonal Bisectional Curvature and the Real Bisectional Curvature.
\end{abstract}

\hfill

In \cite{BroderSBC, BroderSurvey}, the author initiated a program to study the growing number of curvatures in complex geometry by viewing them as quadratic form-valued functions on the unitary frame bundle. This has led to a number of developments, and primary insight which furnished the progress on the Schwarz lemma in \cite{BroderSBC}. The purpose of the present note is to extend this program to the Quadratic Orthogonal Bisectional Curvature (QOBC from here on): $$\text{QOBC}_{\omega} : \mathcal{F}_X \times \mathbb{R}^n \to \mathbb{R}, \hspace*{1cm} \text{QOBC}_{\omega}(v) : = \frac{1}{| v |^2} \sum_{\alpha,\gamma} R_{\alpha \overline{\alpha} \gamma \overline{\gamma}} (v_{\alpha} - v_{\gamma})^2.$$

We recall that the Mori \cite{Mori}, Siu--Yau \cite{SiuYau} solution to the Frankel conjecture informs us that a compact K\"ahler manifold with a K\"ahler metric $\omega$ of positive bisectional curvature $\text{HBC}_{\omega}>0$ is biholomorphic to complex projective space $\mathbb{P}^n$.  Mok \cite{Mok} generalized this result to the case of non-negative bisectional curvature,  further indicating that a sign on the bisectional curvature is a very restrictive assumption. A number of relaxations of the bisectional curvature have been introduced.  

The restriction of the bisectional curvature to orthogonal pairs of tangent vectors yields the \textit{orthogonal bisectional curvature} $\text{HBC}_{\omega}^{\perp}$.  The positivity of $\text{HBC}_{\omega}^{\perp}$ is certainly weaker, algebraically, than the positivity of the full bisectional curvature $\text{HBC}_{\omega}$.  It was shown by Gu--Zhang \cite{GZ}, however, that the K\"ahler--Ricci flow on a compact K\"ahler manifold with $\text{HBC}_{\omega}^{\perp} \geq 0$ converges to a K\"ahler metric with $\text{HBC}_{\omega} \geq 0$.  The classification of compact K\"ahler manifolds with $\text{HBC}_{\omega} \geq 0$ given by Mok's solution \cite{Mok} of the generalized Frankel conjecture applies, showing that the positivity of the orthogonal bisectional curvature furnishes no new examples; they are all biholomorphic to a product of projective spaces and Hermitian symmetric spaces of rank $\geq 2$.

The QOBC was introduced by Wu--Yau--Zheng \cite{WuYauZheng}, where it was shown that on a compact K\"ahler manifold supporting a K\"ahler metric with non-negative QOBC, every nef class is semi-positive (c.f., \cite{DPS}).  Historically, however, it first appears implicitly in the paper of Bishop--Goldberg \cite{BG} as the Weitzenb\"ock curvature operator (c.f., \cite{ChauTam, ChauTam2,PW1,PW2,PW3,PW4}) acting on $(1,1)$--forms.  In contrast with the orthogonal bisectional curvature, the QOBC is strictly weaker than the bisectional curvature,  with an explicit example constructed in \cite{LiWuZheng}.  The QOBC has been the subject of large interest in recent years (see, e.g.,  \cite{BroderGraph,BroderLA,BT1,BT2,ChauTam, ChauTam2,HT, LiWuZheng,Ni,Niu, Tong}).

The purpose of the present short note is to describe the link between combinatorics, distance geometry, and the QOBC.  These results are dispersed throughout the papers \cite{BroderLA, BroderGraph,BT2},  since the considerations here led to solutions of problems in the fields of combinatorics and distance geometry. The papers \cite{BroderLA, BroderGraph, BT2} are intended for those communities, however, and not differential geometers. The present paper is intended for this original audience.

Let us first give a formal definition of the QOBC:

\subsection*{Definition 1}
Let $(X, \omega)$ be a compact K\"ahler manifold.  The QOBC is the function $$\text{QOBC}_{\omega} : \mathcal{F}_X \times \mathbb{R}^n \to \mathbb{R}, \hspace*{1cm} \text{QOBC}_{\omega}(v) : = \frac{1}{| v |^2} \sum_{\alpha,\gamma} R_{\alpha \overline{\alpha} \gamma \overline{\gamma}} (v_{\alpha} - v_{\gamma})^2,$$  where $\mathcal{F}_X$ denotes the unitary frame bundle, $v = (v_1, ...,v_n) \in \mathbb{R}^n$ is a vector, and $\text{QOBC}_{\omega}(0) : = 0$. \\

To state the first theorem, let us introduce some terminology.  A symmetric matrix $A \in \mathbb{R}^{n \times n}$ is declared to be a \textit{Euclidean distance matrix} (EDM) (of embedding dimension $1$) if there is a vector $v = (v_1, ..., v_n) \in \mathbb{R}^n$ such that $A_{\alpha \gamma} = (v_{\alpha} - v_{\gamma})^2$. The well-known Schoenberg criterion \cite{Schoenberg} asserts that a real symmetric matrix with no non-zero diagonal entries (i.e., a hollow matrix) is an EDM if and only if it is negative semi-definite on the hyperplane $H = \{ x \in \mathbb{R}^n : x^t \textbf{e} =0 \}$, where $\textbf{e} = (1,..., 1)^t$.  Let us remark that there are remnants of the Schoenberg criterion in \cite[Lemma 2.6]{ChauTam}, where the K\"ahler form $\omega$ appears to play the role of the vector $\textbf{e}$. An EDM is, by definition, a non-negative matrix in the sense that each entry of the matrix is a non-negative real number.  In particular, the Perron--Frobenius theorem informs us that the largest eigenvalue of $A$ is non-negative and occurs with eigenvector in the non-negative orthant.  This eigenalue is often referred to as the \textit{Perron root of $A$}.  Let $\delta_1 \geq \delta_2 \geq \cdots \geq \delta_n$ denote the eigenvalues of the EDM $A$. We will assume that $\delta_1 \neq 0$. From the above discussion, we know that $\delta_1 >0$ and $\delta_k \leq 0$ for all $k \geq 2$. Define, therefore, for $2 \leq k \leq n$, the \textit{$k$th Perron weight} to be the ratio $r_k : = - \delta_k/\delta_1$, which is a real number lying in the interval $r_k \in [0,1]$.

\subsection*{Theorem 2}
Let $(X^n, \omega)$ be a compact K\"ahler manifold with $\mathcal{R}$ the matrix with entries $\mathcal{R}_{\alpha \gamma} = R_{\alpha \overline{\alpha} \gamma \overline{\gamma}}$. Let $\lambda_1 \geq \lambda_2 \geq \cdots \geq \lambda_n$ denote the eigenvalues of $\mathcal{R}$ with respect to the frame which minimizes the QOBC. Then $\text{QOBC}_{\omega} \geq 0$ if and only if \begin{eqnarray*}
\lambda_1 & \geq & r_2 \lambda_2 + \cdots + r_n \lambda_n
\end{eqnarray*}

holds for all Perron weights $0 \leq r_k \leq 1$. 
\begin{proof}
Fix a frame which minimizes the quadratic orthogonal bisectional curvature of $\omega$.  Let $\lambda_1 \geq \lambda_2 \geq \cdots \geq \lambda_n$ denote the eigenvalues of $\mathcal{R} \in \mathbb{R}^{n \times n}$ and denote by $\delta_1 \geq \delta_2 \geq \cdots \geq \delta_n$ the eigenvalues of an EDM $\Delta$. Write $\mathcal{R} = U^t \text{diag}(\lambda) U$ and $\Delta = V^t \text{diag}(\delta) V$ for the eigenvalue decompositions of $\mathcal{R}$ and $\Delta$. Then  \begin{eqnarray*}
\text{tr}(\mathcal{R} \Delta) \ = \ \text{tr}(U^t \text{diag}(\lambda) U V^t \text{diag}(\delta) V) &=& \text{tr}(VU^t \text{diag}(\lambda) UV^t \text{diag}(\delta) )\\
&=& \text{tr}(Q^t \text{diag}(\lambda) Q \text{diag}(\delta)) \\
&=& \sum_{i,j} \lambda_i \delta_j Q_{ij}^2,
\end{eqnarray*}

where $Q = UV^t$ is orthogonal. The Hadamard square (by which, we mean the matrix $Q \circ Q$ with entries $Q_{ij}^2$) of an orthogonal matrix is doubly stochastic (see, e.g., \cite{Heinz}). The class of $n \times n$ doubly stochastic matrices forms a convex polytope -- the Birkhoff polytope $\mathcal{B}^n$. The minimum of $\text{tr}(\mathcal{R} \Delta)$ is given by \begin{eqnarray*}
\min_{S \in \mathcal{B}^n} \sum_{i,j=1}^n \lambda_i \delta_j S_{ij}.
\end{eqnarray*}
This function is linear in $S$,  achieving its minimum on the boundary of $\mathcal{B}^n$. The well-known Birkhoff--von Neumann theorem tells us that $\mathcal{B}^n$ is the convex hull of the set of permutation matrices, and moreover, the vertices of $\mathcal{B}^n$ are precisely the permutation matrices. Hence,\begin{eqnarray*}
\min_{S \in \mathcal{B}^n} \sum_{i,j=1}^n \lambda_i \delta_j S_{ij} &=& \min_{\sigma \in S_n} \sum_{i=1}^n \lambda_i \delta_{\sigma(i)},
\end{eqnarray*}

where $S_n$ denotes the symmetric group on $n$ letters. An elementary argument (by induction, for instance) shows that \begin{eqnarray*}
\min_{\sigma \in S_n} \sum_{i=1}^n \lambda_i \delta_{\sigma(i)} &=& \sum_{i=1}^n \lambda_i \delta_i.
\end{eqnarray*}

From the discussion after Definition 1, this completes the proof.
\end{proof}

\subsection*{Remark 3}
The main theorem is interesting for a number of reasons. The first is that there is an eigenvalue characterization in terms of the matrix $\mathcal{R}$. Of course, this matrix requires a frame to be fixed, but given that there are a number of frame-dependent curvatures in complex geometry that have appeared in recent years (most notably, the real bisectional curvature \cite{YangZhengRBC}, and the Schwarz bisectional curvatures \cite{BroderSBC, BroderSurvey}), this offers some insight into the relationship between these curvatures.  Further detail in this direction will be discussed shortly. \\

The existence of an eigenvalue characterization is surprising in itself \cite{BroderLA} since Euclidean distance matrices are defined in a frame-dependent manner; the class of positive matrices (i.e., matrices with positive entries) are certainly not invariant under a change of basis.  We suspect that this result (motivated entirely by complex-geometric considerations) will have further generalizations and applications to combinatorics and distance geometry. \\

The class of matrices $A$ satisfying $\sum_{\alpha,\gamma=1}^n A_{\alpha \gamma} (v_{\alpha} - v_{\gamma})^2 \geq 0$ form the so-called \textit{dual EDM cone}.  This is a completely elementary observation: An EDM (of embedding dimension $1$) is a matrix of the form $B_{\alpha \gamma} :=  (v_{\alpha} - v_{\gamma})^2$.  Hence, we can write $\sum_{\alpha,\gamma=1}^n A_{\alpha \gamma} (v_{\alpha} - v_{\gamma})^2 \geq 0$ as $\text{tr}(AB) \geq 0$, from which we immediately see the following: 

\subsection*{Proposition 4}
Let $(X, \omega)$ be a compact K\"ahler manifold.  Then $\text{QOBC}_{\omega} \geq 0$ if and only if, with respect to the frame which minimizes the QOBC, the matrix $\mathcal{R}$ lies in the dual EDM cone. \\

The above result, albeit elementary, is important, in that it gives us the appropriate language to speak when considering the QOBC.  It also allows us to exploit the results of distance geometry and combinatorics to say something about the QOBC. For instance, using Dattorro's dual EDM cone criterion \cite{D1, D2}, we have:

\subsection*{Theorem 5}\label{Remark}
Let $\delta : \mathbb{R}^n \to \mathcal{S}_{\text{diag}}^n$ be the operator mapping a vector $v \in \mathbb{R}^n$ to the diagonal matrix $\text{diag}(v)$. Then a real symmetric matrix $A$ lies in the dual EDM cone if and only if $\delta(A \textbf{e}) - A$ is positive-semi-definite. In particular,  $$\text{QOBC}_{\omega} \geq 0 \ \iff \ \delta(\mathcal{R} \textbf{e}) - \mathcal{R} \in \mathcal{PSD}.$$ 

\subsection*{Remark 6}

Recall that the real bisectional curvature $\text{RBC}_{\omega}$ of a Hermitian metric $\omega$ was introduced by Yang--Zheng \cite{YangZhengRBC} as follows: \begin{eqnarray*}
\text{RBC}_{\omega} : \mathcal{F}_X \times \mathbb{R}^n \to \mathbb{R}, \hspace{1cm} \text{RBC}_{\omega}(v) : = \sum_{\alpha, \gamma=1}^n R_{\alpha \overline{\alpha} \gamma \overline{\gamma}} v_{\alpha} v_{\gamma}.
\end{eqnarray*}

The real bisectional curvature is clearly similar in appearance to the QOBC.  The above theorem, however, indicates that they are in fact \textit{opposite} in the following sense: The condition that $\text{RBC}_{\omega} > 0$ translates to the matrix $\mathcal{R}$ being positive-definite. Let $\mathcal{PD}$  denote the cone of positive-definite matrices, and let $\mathcal{EDM}$ denote the cone of Euclidean distance matrices. The $\mathcal{PD}$ cone is self-dual, while the $\mathcal{EDM}$ cone is not self-dual, and moreover, from \nameref{Remark}, the intersection of these cones is trivial: $$\mathcal{PD} \cap \mathcal{EDM} = \{ 0 \}.$$

Hence, the real bisectional curvature (corresponding to $\mathcal{R} \in \mathcal{PD}$) is \textit{opposite} to the QOBC (corresponding to $\mathcal{R} \in \mathcal{EDM}^{\ast})$ in the sense that the dual cones intersect trivially.

\subsection*{Combinatorics and Distance Geometry}
Let us conclude by describing the curious link with combinatorics and distance geometry.  Indeed,  let $G$ be a finite weighted graph (with possibly negative weights) with vertex set $V(G) = \{ v_1, ..., v_n \}$.  Let $A \in \mathbb{R}^{n \times n}$ be the adjacency matrix specifying the weighting. The Dirichlet energy for a weighted graph $(G,A)$ is defined by \begin{eqnarray*}
\mathcal{E}(f) & : = & \sum_{\alpha, \gamma=1}^n A_{\alpha \gamma} (f(v_{\alpha}) - f(v_{\gamma}))^2,
\end{eqnarray*}

where $f : V(G) \to \mathbb{R}$ is a function defined on the vertices of $G$.

\subsection*{Theorem 7}
Let $(X, \omega)$ be a compact K\"ahler manifold.  Then $\text{QOBC}_{\omega} \geq 0$ is equivalent to the non-negativity of the Dirichlet energy of every weighted graph $(G,A)$ with $G$ a finite graph with $n$ vertices and $A$ given by the matrix $\mathcal{R} = (\mathcal{R}_{\alpha \gamma})$. 

\subsection*{Remark 8}
Recall that the QOBC first appears in the paper of Bishop--Goldberg \cite{BG} as the Weitzenb\"ock curvature operator (c.f., \cite{PW1, PW2, PW3, PW4}) acting on real $(1,1)$--forms.  In other words, if $\Delta_g$ denotes the Bochner Laplace operator, and $\Delta_d$ denote the Laplace--Beltrami operator, acting on real $(1,1)$--forms,  their difference realizes the QOBC.  What is curious, is that the above theorem indicates that the difference of these Laplace operators is the (discrete analog of the) Dirichlet energy associated with the curvature. We hope that those more experienced in the discrete theory can give further insight in this direction.

\subsection*{Acknowledgements}
The author would like to thank his advisors Ben Andrews and Gang Tian for their support, encouragement, and interest in this work.  He would also like to express his gratitude to Peter Petersen, Fangyang Zheng,  Yanir Rubinstein,  Renato Bettiol, and Kai Tang for their interest and valuable communications.

\subsection*{Data availability}
Data sharing not applicable to this article as no datasets were generated or analysed during the current study.

\subsection*{Declarations}

\subsection*{Conflict of interest}
On behalf of all authors, the corresponding author states that there is no conflict of interest.

\end{document}